\documentclass[11pt]{article}
\usepackage{latexsym}

\newcommand{\fieldk}{\hbox{{\rm k}}}

\newcommand{\link}{\hbox{\rm link}}
\newcommand{\im}{\hbox{\rm Im}}

\newcommand{\WDVV}{\hbox{$W$}}
\newcommand{\preWDVV}{\hbox{$R$}}

\newcommand{\Hilb}{{\mathcal H}}

\newcommand{\qed}{\mbox{$\Box$}\vspace{\baselineskip}}

\newenvironment{proof}{\noindent {\bf Proof:}}{{\qed}}
\newenvironment{proof_special}{\noindent {\bf Proof:}}{}

\newtheorem{theorem}{Theorem}[section]
\newtheorem{proposition}[theorem]{Proposition}
\newtheorem{lemma}[theorem]{Lemma}

\newtheorem{corollary}[theorem]{Corollary}

\begin{document}

\title{The pre-WDVV ring of physics and its topology
\footnote{
To appear in 
{\em The Ramanujan Journal}, {\bf 10},
No.\ 2 (2005),
269--281.
2000 Mathematics Subject Classification:  Primary:  13F55,
Secondary:  05E99 and 55P15.
Keywords:  WDVV equations, 
Keel's presentation, 
moduli space,
pre-WDVV complex,
Morse matching,
Cohen-Macaulay,
Whitehouse complex.}
}
\author{{\sc Margaret A.\ READDY}\\\\
{\small Dedicated to Louis Billera in honor of his 60th Birthday}
}

\date{}

\maketitle

\begin{abstract}
We show how 
 a simplicial complex arising from
the WDVV (Witten-Dijkgraaf-Verlinde-Verlinde)
equations of string theory
is the Whitehouse complex.
Using discrete Morse theory,
we give an elementary proof that
the Whitehouse 
complex $\Delta_n$ is
homotopy equivalent to a wedge of $(n-2)!$ spheres
of dimension $n-4$.
We also verify the Cohen-Macaulay property.
Additionally,
recurrences are given for the face enumeration
of the complex and the Hilbert series
of the associated pre-WDVV ring.
\end{abstract}

\section{Introduction}

The moduli space of smooth $n$-pointed stable curves
of genus $g$, denoted $\overline{M}_{g,n}$,
was introduced by 
Deligne, Mumford and 
Knudsen~\cite{Deligne_Mumford, Knudsen, Mumford_stability}
to give a natural compactification
of Mumford's~\cite{Mumford}  moduli space of nonsingular
curves of genus $g$.
A new construction of the
genus $0$ case
is due to 
Keel~\cite{Keel} 
using blowups.
Keel also gives
the presentation
for the generators and relations of
the cohomology ring of $\overline{M}_{0,n}$.

The Associativity Equations in physics,
also known as the WDVV Equations~\cite[Section 0.2]{Kontsevich_Manin_II},
are a system of partial differential equations due
to Witten, Dijkgraaf, Verlinde and 
Verlinde~\cite{Dijkgraaf_Verlinde_Verlinde, Dubrovin, Witten}.
Their solutions are generating functions
$\Phi$ 
having coefficients which encode ``potential'' or ``free energy''
and
are related to understanding 
quantum gravity
via the topological quantum theory approach~\cite{Witten_I}.

In~\cite{Kontsevich_Manin_I}
Kontsevich and Manin determine
the potential function~$\Phi$ in the case of Fano manifolds.
Underlying their work is the need to construct Gromov-Witten
classes.  These are
linear maps between the cohomology
of a projective algebraic
manifold and 
the cohomology of the moduli space~$\overline{M}_{g,n}$.
Kontsevich and Manin develop a cohomological field theory 
in terms of 
the Gromov-Witten class language and the language of 
operads.  From this theory
and Keel's presentation,
they derive all the
linear relations 
between homology classes of boundary strata of
any codimension~\cite[Sections~6 and 7]{Kontsevich_Manin_I}.
Since Keel's presentation and 
the splitting axiom for the
Gromov-Witten 
classes~\cite[equations (0.3) and (0.4)]{Kontsevich_Manin_II})
imply the WDVV equations (see~\cite{Kontsevich_Manin_I} for details),
it is natural to refer to the cohomology ring given by
Keel's presentation 
as the WDVV ring.

In a previous paper, the author studied an analogue of
the WDVV ring~\cite{Readdy}.  This ring,
known as the Losev-Manin ring, arises instead from 
the Commutativity Equations
of physics~\cite{Losev_Manin}.
Like the WDVV ring, the defining ideal for
the Losev-Manin ring is composed of linear and quadratic relations.
Since the quadratic relations determine
the interesting behavior of 
the Losev-Manin ring,
it is again natural to take a closer look at the quadratic relations
defining the WDVV ring.
We call this new ring the pre-WDVV ring.

In this paper we show  the pre-WDVV ring
can be realized as
the Stanley-Reisner ring
of the Whitehouse complex
$\Delta_n$.
This is a well-studied complex 
which goes back to work of Boardman~\cite{Boardman}.
See Section~\ref{section_definitions} for further references.
As a result, many of our results are rediscoveries, but
with simpler proofs.
Since our overall goal is
to find a natural combinatorial
object corresponding to the minimal generators  of the WDVV ring,
we expect our current approach to 
better lead to its understanding.

In Section~\ref{section_pre_WDVV} we study the links of faces
in the Whitehouse complex  via a forest representation of the faces.
We then describe five injective maps which map
$\Delta_n$ onto $\Delta_{n+1}$.
As new results, this allows us to deduce recurrences for the face vectors 
and Hilbert series of $\Delta_n$.
In Section~\ref{section_Morse_matching}
we construct a Morse matching
on the complex
to give an elementary proof that
the Whitehouse complex
$\Delta_n$ is homotopy equivalent to
a wedge of $(n-2)!$ spheres 
of dimension $n-4$; see Theorem~\ref{theorem_pre-WDVV_homotopy}.
In Section~\ref{section_Cohen_Macaulayness}
we again return to the geometric description
of the links of faces in $\Delta_n$ 
to give an elementary proof that the Whitehouse 
complex is Cohen-Macaulay.
In the last section we indicate some work in progress about
the related WDVV ring.

\section{Definitions and notation}
\label{section_definitions}

Throughout we will assume familiarity with
basic poset terminology and combinatorial concepts.
An excellent reference for the uninitiated is Stanley's
text~\cite{Stanley}.

Let
$V$ be a finite vertex set.
A simplicial complex $\Delta$
is
a collection of subsets of $V$ called faces 
satisfying
$\emptyset \in V$ 
and
$\{v\} \in \Delta$ for all $v \in V$,
and
if
$F \subseteq G \in \Delta$ then $F \in \Delta$.
The {\em dimension} of a face $F \in \Delta$
is given by $|F| - 1$.
A face $F$ in $\Delta$ 
is a {\em facet} if there is no  face
in $\Delta$ strictly containing $F$,
while a simplicial complex is 
{\em pure}
if all the facets have the same dimension.
The {\em link} of
a face $F$ in the complex $\Delta$ is defined to be
$\link_{\Delta}(F) = \{ G \subseteq V \:\: : \:\: F \cap G = \emptyset
\mbox{     and     }
F \cup G \in \Delta\}$.
Finally,
for $\Delta$ and $\Gamma$ two simplicial complexes having disjoint vertex sets
$V$ and $W$, the
{\em join} of $\Delta$ and $\Gamma$,
denoted
$\Delta * \Gamma$,
is the simplicial complex
on the vertex set
$V \cup W$ consisting of the faces
$\Delta * \Gamma = \{ F \cup G \:\: : \:\: F \in \Delta 
\mbox{  and  } G \in \Gamma\}$.

For completeness, 
we give the physicist's definition
of  the WDVV ring~\cite[Section 0.10]{Kontsevich_Manin_II}
based on Keel's presentation~\cite[Theorem 1]{Keel}.
Let $n \geq 3$
and let $\fieldk$ be a field of characteristic zero.
We say $\sigma$ is a {\em stable $2$-partition of
$\{1, \ldots, n\}$} if
$\sigma$ is an unordered partition of the elements
$\{1, \ldots, n\}$ into two blocks
$\sigma = S_1 / S_2$ with
$|S_i| \geq 2$.
Denote by $P_n$ the set of stable $2$-partitions
of $\{1, \ldots, n\}$.
For each $\sigma \in P_n$,
the element $x_{\sigma}$ corresponds to
a cohomology class of 
$H^*(\overline{M}_{0,n})$,
that is,
the cohomology ring of 
the moduli space of
$n$-pointed stable curves of genus zero.
Throughout it will be convenient for fixed~$\sigma$ 
to think of $x_{\sigma}$ as simply an indeterminate.
For $\sigma = S_1/ S_2$ and
$\tau =  T_1/ T_2$ from the index set $P_n$,
let
$a(\sigma,\tau)$
 be  the number of nonempty pairwise
distinct sets among $S_i \cap T_j$,
where $1 \leq i, j \leq 2$.
Define the ideal~$I_n$ in the
polynomial ring
$\fieldk[x_{\sigma} \: : \: \sigma \in P_n]$
by:
\begin{enumerate}
\item
(linear relations)
For $i,j,k,l$ distinct:
\begin{equation}
\label{equation_linear_relations}
        R_{ijkl} = \sum_{ij \sigma kl} x_{\sigma}
                    -
                    \sum_{kj \tau il} x_{\tau},
\end{equation}
where the summand
$ij \sigma kl$ means to sum over all stable $2$-compositions
$\sigma = S_1 /S_2$ with
the elements
$i, j \in S_1$ and the elements
$k, l \in S_2$.

\item
(quadratic relations)
For each pair $\sigma$ and $\tau$ with $a(\sigma,\tau) = 4$,
\begin{equation}
\label{equation_quadratic_relations}
        x_{\sigma} \cdot x_{\tau}
\end{equation}
\end{enumerate}
The 
{\em WDVV ring}
(or {\em Associativity ring})
$\WDVV_n$ is defined to be
the quotient ring
$\WDVV_n = \fieldk[x_{\sigma} \: : \: \sigma \in P_n]/I_n$.

Rather than working with the complex associated with the
the entire WDVV ring, we now consider the 
pre-Associativity or pre-WDVV ring defined by
taking the WDVV ring modulo
the quadratic relations only.
More formally,
the {\em pre-WDVV ring}
$\preWDVV_n$
is defined to be
$\preWDVV_n = \fieldk[x_{\sigma} \: : \: \sigma \in P_n]/{J_n}$,
where $J_n$ is the ideal generated by the quadratic
relations~(\ref{equation_quadratic_relations}) of $I_n$.

Keel's presentation of the cohomology ring
of the moduli space 
$\overline{M}_{0,n}$ has twice as many variables
as he instead indexes
the cohomology classes by 
subsets $S$ of $\{1, \ldots, n\}$.
However, Keel makes the further requirement that
$x_S = x_{\overline{S}}$,
where
$\overline{S}$ is the complement of
$S$ taken with respect to $\{1, \ldots, n\}$,
making his variables correspond to
stable $2$-partitions.
Hence, it is natural to think of the
variables $x_{\sigma}$ for
$\sigma \in P_n$ to be indexed instead by
subsets
$S \subseteq \{2, \ldots, n\}$
having cardinality
satisfying $2 \leq |S|\leq n-2$.
See~\cite[Theorem 1]{Keel}.

Based upon Keel's presentation,
we next 
construct a simplicial complex 
intimately related to the pre-WDVV ring.
Fix an integer $n \geq 3$.
Let
$\Delta_n$ 
be the simplicial complex
with vertex set
$V_n = \{S \subseteq \{2, \ldots, n\} 
\:\: : \: \: 2 \leq |S| \leq n-2 \}$.
A subset $F \subseteq V_n$ is a face if
for all $S, T \in F$ either
$S \subseteq T$, $S \supseteq T$ or $S \cap T = \emptyset$.

Observe that the complex $\Delta_n$ 
is described by its minimal non-faces, and that
each minimal non-face has cardinality $2$.
Also,
a set $S$ in the vertex set
of $\Delta_n$ corresponds
to a stable $2$-partition
$S / \overline{S}$,
where the  complement is taken with respect
to the set $\{1, \ldots, n\}$.
Hence the condition of being a face 
$F$ corresponds to
$a(\sigma,\tau) < 4$ 
for
all $S, T \in F$ where
$\sigma = S / \overline{S}$
and
$\tau = T / \overline{T}$.

The complex $\Delta_n$
coincides 
with 
Boardman's space of fully grown $n$-trees~\cite{Boardman}.
This space was rediscovered by Whitehouse in her 
thesis~\cite{Whitehouse}
and independently by Culler and Vogtmann~\cite{Culler_Vogtmann}.
In the literature it is commonly referred to as the Whitehouse complex.
The Whitehouse complex also has connections with phylogenetic trees.
See the work of Billera, Holmes and Vogtmann~\cite{Billera_Holmes_Vogtmann}.

Recall for a finite simplicial complex $\Delta$ 
with vertex set $V$ the
{\em Stanley-Reisner ring} $\fieldk[\Delta]$ is
defined to be 
$\fieldk[\Delta] = \fieldk[x_v \: : \: v \in V]/I(\Delta)$,
where $I(\Delta)$
is the ideal generated 
by the non-faces of 
the complex $\Delta$; see~\cite{Stanley_green}.
Using this terminology,
we immediately have the following result.

\begin{proposition}
The pre-WDVV ring $\preWDVV_n$ is the
Stanley-Reisner ring of the Whitehouse complex~$\Delta_n$,
that is,
$\preWDVV_n = \fieldk[\Delta_n]$.
\end{proposition}

%
%

\newcommand{\ellipse}[6]{
\qbezier(#1,#5)(#1,#4),(#2,#4)
\qbezier(#2,#4)(#3,#4)(#3,#5)
\qbezier(#3,#5)(#3,#6)(#2,#6)
\qbezier(#2,#6)(#1,#6)(#1,#5)
}

\begin{figure}[t]
\begin{picture}(400,120)


\ellipse{40}{70}{100}{30}{45}{60}
\ellipse{110}{170}{230}{40}{65}{90}
\ellipse{125}{145}{165}{50}{65}{80}

\put(58,41){2}
\put(68,41){3}
\put(78,41){4}

\put(137,63){5}
\put(147,63){6}

\put(187,63){7}

\put(68,81){8}


\put(300,110){\line(-1,-2){20}}
\put(300,110){\line(0,-2){40}}
\put(300,110){\line(1,-2){20}}

\put(277,55){2}
\put(298,55){3}
\put(318,55){4}

\put(300,110){\circle*{4}}
\put(300,70){\circle*{4}}
\put(280,70){\circle*{4}}
\put(320,70){\circle*{4}}

\put(370,70){\line(-1,-2){20}}
\put(370,70){\line(1,-2){20}}

\put(370,70){\circle*{4}}
\put(350,30){\circle*{4}}
\put(390,30){\circle*{4}}

\put(348,15){5}
\put(388,15){6}

\put(390,110){\line(-1,-2){20}}
\put(390,110){\line(1,-2){20}}

\put(408,55){7}

\put(390,110){\circle*{4}}
\put(410,70){\circle*{4}}

\put(410,70){\circle*{4}}

\put(440,110){\circle*{4}}

\put(438,95){8}

\end{picture}
\caption{The forest representation of the face  
$F = \{\{2,3,4\}, \{5,6\}, \{5,6,7\} \}$
in
$\Delta_8$.}
\label{figure_forest_example}
\end{figure}
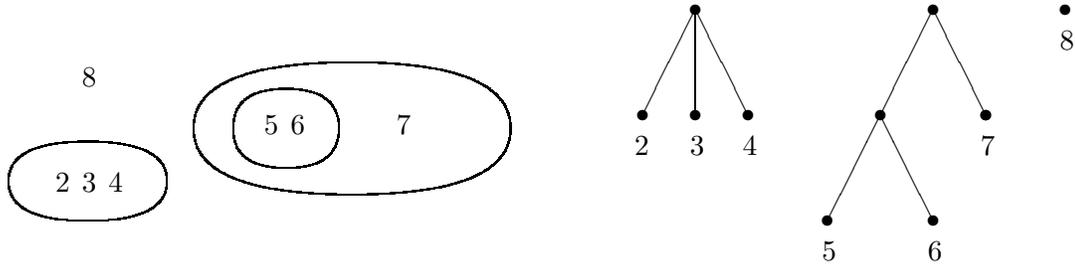

\section{Facial structure}
\label{section_pre_WDVV}

We now proceed with a more formal study of the
Whitehouse complex $\Delta_n$
which will lead to understanding its topology.
In particular, as new results we derive recurrences for its
face $f$-vector and $h$-vector.

Observe the complex
$\Delta_3$ consists of the empty set.
It will be convenient to view this complex as a
$(-1)$-dimensional sphere.
The complex $\Delta_4$ consists of $3$ isolated vertices,
while the complex~$\Delta_5$ 
is the Peterson graph.
Note the two last examples are
homotopy equivalent to the wedge of two $0$-dimensional
spheres,
respectively,
the wedge of six $1$-dimensional spheres.

%
%

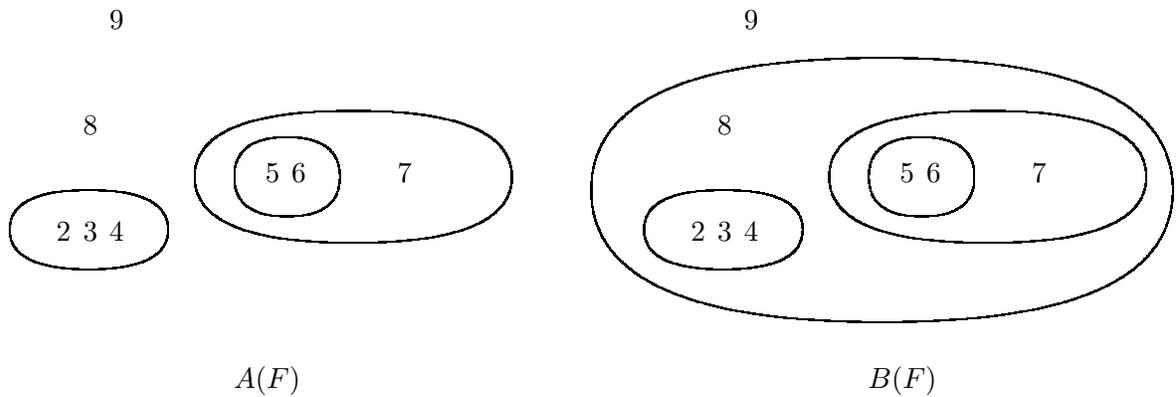
\begin{figure}[t]
\begin{picture}(460,200)



\ellipse{20}{50}{80}{70}{85}{100}
\ellipse{90}{150}{210}{80}{105}{130}
\ellipse{105}{125}{145}{90}{105}{120}

\put(38,81){2}
\put(48,81){3}
\put(58,81){4}

\put(117,103){5}
\put(127,103){6}

\put(167,103){7}

\put(48,121){8}
\put(58,161){9}

\put(105,25){$A(F)$}


\ellipse{240}{350}{460}{50}{100}{150}
\ellipse{260}{290}{320}{70}{85}{100}
\ellipse{330}{390}{450}{80}{105}{130}
\ellipse{345}{365}{385}{90}{105}{120}

\put(278,81){2}
\put(288,81){3}
\put(298,81){4}

\put(357,103){5}
\put(367,103){6}

\put(407,103){7}

\put(288,121){8}
\put(298,161){9}

\put(345,25){$B(F)$}

\end{picture}
\caption{The maps $A$ and $B$
applied to the face 
$F = \{ \{2,3,4\}, \{5,6\}, \{5,6,7\} \}$
in $\Delta_8$.}
\label{figure_maps_A_B}
\end{figure}

Each face $F$ of $\Delta_n$ can be described by a forest.
The leaves of the forest
are $2, \ldots, n$.  
The internal nodes are the elements of the face $F$.
The cover relation of the forests
is defined as follows:
For $S, T \in F$, we say
$S$ covers $T$ if $S \supset T$ 
and there is no $U \in F$ such that
$S \supset U \supset T$.
Moreover,
$S \in F$ covers
a leaf
$i \in \{2, \ldots, n\}$
if there is no
$U \in F$ such that
$i \in U \subset S$.
(Here we use $\subset$ and
$\supset$ to mean strict 
containment and reverse containment, respectively.)
See Figure~\ref{figure_forest_example}
for the example
$F = \{ \{2, 3, 4\}, \{5, 6\}, \{5, 6, 7\} \}$. From 
the forest representation
it is straightforward to see that the
complex $\Delta_n$ is pure of dimension $n-4$.

We have the following results about the links of faces
in the complex $\Delta_n$.

\begin{lemma}
\label{lemma_vertex_link}
Let $T$ be a vertex of the Whitehouse complex $\Delta_n$.
Then
$$
        \link_{\Delta_n}(T) \cong \Delta_{|T| + 1} * \Delta_{n- |T| + 1},
$$
where
$*$ denotes the join operation.
\end{lemma}
\begin{proof}
First observe that
for $T$ a vertex of $\Delta_n$,
the link of $T$ 
consists of all faces that contain~$T$.
Each face in
the $\link_{\Delta_n}(T)$
can be described as a forest on the elements
$\{2, \ldots, n\}$.
In a given forest, 
the element $T$
is an internal node.
One can choose the structure beneath $T$ in the forest
in $\Delta_{|T| + 1}$ ways.
To build the rest of the forest,
treat the tree built by $T$ as a leaf.
Then
one sees 
there are
$n - 1 - |T| + 1$ elements
available.
The resulting structure has
the form
$\Delta_{n - |T| + 1}$.
Hence,
$\link_{\Delta_n}(T) \cong \Delta_{|T| + 1} * \Delta_{n- |T| + 1}$.
\end{proof}

%
%

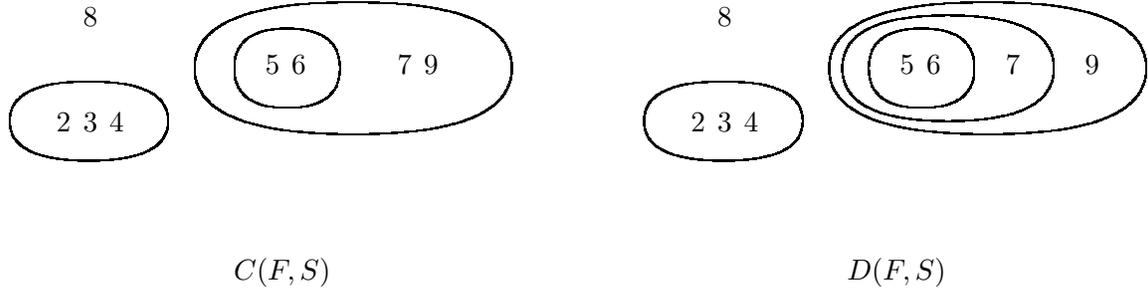
\begin{figure}[t]
\begin{picture}(400,180)



\ellipse{20}{50}{80}{70}{85}{100}
\ellipse{90}{150}{210}{80}{105}{130}
\ellipse{105}{125}{145}{90}{105}{120}

\put(38,81){2}
\put(48,81){3}
\put(58,81){4}

\put(117,103){5}
\put(127,103){6}

\put(167,103){7}
\put(177,103){9}

\put(48,121){8}

\put(105,25){$C(F,S)$}


\ellipse{260}{290}{320}{70}{85}{100}
\ellipse{330}{390}{450}{80}{105}{130}

\ellipse{335}{375}{415}{85}{105}{125}

\ellipse{345}{365}{385}{90}{105}{120}

\put(278,81){2}
\put(288,81){3}
\put(298,81){4}

\put(357,103){5}
\put(367,103){6}

\put(397,103){7}
\put(427,103){9}

\put(288,121){8}

\put(337,25){$D(F,S)$}

\end{picture}

\caption{The maps $C$ and $D$
applied to the face 
$F = \{ \{2,3,4\}, \{5,6\}, \{5,6,7\} \}$
and
set $S = \{5,6,7\}$
in $\Delta_8$.}
\label{figure_maps_C_D}
\end{figure}

Iterating 
Lemma~\ref{lemma_vertex_link}
gives the following more general result.

\begin{proposition}
\label{proposition_link}
Let $F$ be a face of the Whitehouse complex $\Delta_n$.
Then
$$
        \link_{\Delta_n}(F) \cong \Delta_{c + 1} * 
                            \prod_{T \in F}
                            \Delta_{c(T) + 1},
$$
where
$\prod$ denotes taking the join operation $*$ among the 
factors, $c$ is the number of components in the 
forest representation of $F$,
and 
$c(T)$ is the number of children the node $T$ has in
the forest representation of $F$.
\end{proposition}

We introduce five maps
$A, B, C, D, E \: \: : \: \: \Delta_n \rightarrow \Delta_{n+1}$
which together map
$\Delta_n$ onto $\Delta_{n+1}$.
For a face $F \in \Delta_n$ define
$$A(F) = F
\:\:\:
\mbox{        and        }
\:\:\:
B(F) = F \cup \{ \{2, \ldots, n\} \}.
$$
For $S \in F$ define
$$C(F,S) = \{ T \cup \{ n+1\} \: \: : S \subseteq T, T \in F\}
        \cup
        \{T \: \: : \: \: S \not\subseteq T, T \in F\}
$$
and
$$
        D(F,S) = C(F,S) \cup \{ S \}.
$$
Finally, for $2 \leq i \leq n$ define
$$
        E(F,i) =
        \{ \{ i, n+1\}\} \cup
        \{ T \cup \{n+1\} \: \: : \: \: i \in T, T \in F \}
        \cup
        \{T \: \: : \: \: i \not\in T, T \in F\}.
$$
Informally speaking,
the $A$ map adds the element $n+1$ to 
the set diagram of the face $F$
whereas 
the $B$ map creates
a new set 
consisting of the entire set diagram of the face $F$ and
then adds the element $n+1$.
The $C$ map selects a set $S$ from the face $F$
and adds the element
$n+1$ to it,
whereas the $D$ map selects a set $S$ from the face $F$
and creates a new set consisting of the set
$S$ and the element $n+1$.
The $E$ map replaces an element $i$ with the
set consisting of $i$ and $n+1$.
For $F$ a face of dimension $i$
the maps $A$ and $C$ leave the dimension
of $F$ unchanged while the maps
$B$, $D$ and $E$ each increase the dimension
by one.

\begin{figure}[t]
\begin{picture}(255,160)


\ellipse{130}{160}{190}{70}{85}{100}
\ellipse{200}{260}{320}{80}{105}{130}
\ellipse{215}{235}{255}{90}{105}{120}
\ellipse{265}{285}{305}{90}{105}{120}

\put(148,81){2}
\put(158,81){3}
\put(168,81){4}

\put(227,103){5}
\put(237,103){6}

\put(277,103){7}
\put(287,103){9}

\put(158,121){8}

\put(217,25){$E(F,7)$}

\end{picture}
\caption{The map $E(F,7)$
where 
$F = \{ \{2,3,4\}, \{5,6\}, \{5,6,7\} \}$
in $\Delta_8$.}
\label{figure_map_E}
\end{figure}
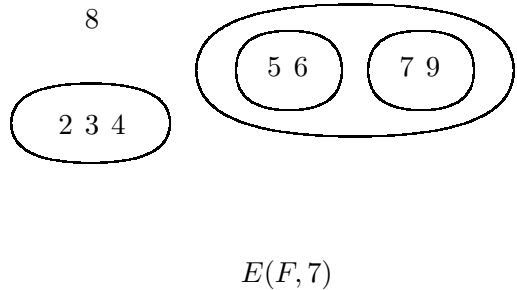

It is easily seen that all of these maps are injective and that
$\Delta_{n+1}$ is a disjoint union of their images.
See Figures~\ref{figure_maps_A_B},
\ref{figure_maps_C_D} and \ref{figure_map_E}
for an example of each of these maps.

Observe that for fixed $i$ the image
of $E(F,i)$ 
as $F$ runs over all faces
$F$ in the complex
$\Delta_n$
is the link of the vertex
$\{i, n+1\}$ 
in $\Delta_{n+1}$ which
is isomorphic to $\Delta_n$.
Hence we can decompose the face
poset of the complex
$\Delta_{n+1}$ into the
images of the five maps $A$, $B$, $C$, $D$ and 
$E$.
In turn, the image of $E$ 
further decomposes into $n-1$ copies
of the face poset of $\Delta_n$.  
It is straightforward to see each 
of the $n-1$ copies  of 
the $\Delta_n$ in the face
poset of~$\Delta_{n+1}$ is an upper order ideal.

%
%

 From the maps we have just defined,
we can give a recurrence for the number of faces in the 
Whitehouse complex.

\begin{theorem}
\label{theorem_face_recursion}
Let $f_{n,i}$ denote the number of faces in
the Whitehouse complex $\Delta_n$ having dimension~$i$.
Then
$f_{n,i}$ satisfies the recursion
$$
        f_{n,i} = (i+2) \cdot f_{n-1,i} + (n+i-1) \cdot f_{n-1,i-1},
$$ 
where $-1 \leq i \leq n-4$
and
$f_{n,-1} = 1$.
\end{theorem}
\begin{proof}
Easily $f_{n,-1} = 1$ for all $n \geq 3$, as there is one
$(-1)$-dimensional face in $\Delta_n$, namely
the empty set.
To show the recurrence,
note that a face having dimension $i$ has
$i+1$ sets occurring in its set representation.
The maps $A, B, C, D$ and $E$ each
either
leave the number of sets unchanged or increase the number
by one.
Since they map $\Delta_{n-1}$ bijectively onto $\Delta_n$,
we can build 
all the $i$-dimensional faces in $\Delta_n$.
First, we take a face from $\Delta_{n-1}$ 
that has $i+1$ sets.
The $A$ map simply adds the element $n$ in $1$ way,
while the $C$ map adds it in the number of sets ways,
that is, $i+1$ ways .  Overall we have constructed
$(i+2) \cdot f_{n-1,i}$ new faces in $\Delta_n$.
The other way to build an $i$-dimensional face
is to take a face from 
$\Delta_{n-1}$ having $i$ sets and add the element
$n$ in such a way to increase the number of sets.
The $B$ map does this in $1$ way,
the $D$ map in number of sets ways, that is, $i$ ways,
and the
$E$ map does this in $n-2$ ways.
Hence we have constructed
$(n+i-1) \cdot f_{n-1,i-1}$ new faces
in $\Delta_n$,
and the recurrence holds.
\end{proof}

\begin{corollary}
\label{corollary_facets}
The Whitehouse complex
$\Delta_n$ is
pure of dimension $n-4$ with
$(2n-5)!! = (2n-5)\cdot (2n-7) \cdots 3 \cdot 1$
facets.
\end{corollary}
This corollary can be found in~\cite{Robinson_Whitehouse}
and~\cite{Vogtmann}.

Recall the {\em $h$-vector}
of a $(d-1)$-dimensional simplicial complex
$\Delta$
is the vector
$(h_{0}, \ldots, h_{d})$
where
$h_{k} = \sum_{i=0}^k (-1)^{k-i} {d-i \choose k-i} \cdot f_{d-1,i-1}$.
The Hilbert series of the Stanley-Reisner ring of
$\Delta$ is thus given by
${\cal H}(\fieldk[\Delta])
 = (h_0 + h_1 \cdot t + \cdots + h_{d} \cdot t^d)/(1-t)^d$.
See~\cite{Stanley_green}.
Using the expression for the $h$-vector in terms of the
$f$-vector and Theorem~\ref{theorem_face_recursion},
we obtain
the following identity.

\begin{corollary}
\label{corollary_h_vector_recursion}
The $h$-vector
of the Whitehouse complex $\Delta_n$ satisfies
the recurrence
$$
	h_{n,k} = (k+1) \cdot h_{n-1,k} + (2n-k-5) \cdot h_{n-1,k-1},
$$
where $0 \leq k \leq n-3$.
\end{corollary}
The 
Hilbert series of the pre-WDVV ring $\preWDVV_n$ 
for $3 \leq n \leq 8$ (and thus the $h$-vector values for 
the Whitehouse complex $\Delta_n$)
are
displayed in Table~\ref{table_pre_WDVV_ring}.

%
%

\begin{table}
\begin{center}
\begin{tabular}{l|l}
$n$  &  $(1-t)^{n-3} \cdot \Hilb(\preWDVV_n)$\\  \hline
$3$ &   $1$\\
$4$ & $1 + 2t$\\
$5$ & $1 + 8t + 6t^2$\\
$6$ & $1 + 22t + 58t^2 + 24 t^3$\\
$7$ & $1 + 52t + 328 t^2 + 444 t^3 + 120 t^4$\\
$8$ & $1 + 114 t + 1452 t^2 + 4400 t^3 + 3708 t^4 + 720 t^5$
\end{tabular}
\end{center}
\caption{The Hilbert series of the pre-WDVV ring $\preWDVV_n$ for
$3 \leq n \leq 8$.}
\label{table_pre_WDVV_ring}
\end{table}

\section{Morse matching and the topology of $\Delta_n$}
\label{section_Morse_matching}

R.\ Forman
devised a discrete version of Morse theory to
study the topology of simplicial complexes,
and more generally,
CW-complexes.
We give a brief overview of this.
More details can be found in~\cite{Forman}.
We then describe a Morse matching of the Whitehouse complex
and 
verify its homotopy type.

Let~$P$ be an arbitrary poset.
Begin by orienting all of the edges
in the Hasse diagram of $P$,
that is,
the
cover relations of $P$, downwards.  
Next, form a matching $M$ on the elements of $P$.
Reverse the orientation of the edges in the matching
to be upwards.
Such a matching $M$ is a {\em Morse matching}
if the resulting directed graph is acyclic.
An unmatched element
of $P$ is called
{\em critical}.

For a simplicial complex $\Delta$ the {\em face poset}
$P$ is the poset formed by taking the faces of the
complex as elements and ordering them by inclusion.
The face poset is ranked with
the rank of an element $x \in P$
given by $i$ if
$|x| = i$.
In the case we do not wish to include the empty face
in the face poset, we denote the resulting face poset by
$\overline{P}$.

The following is Forman's result~\cite{Forman}.

\begin{theorem}
For a simplicial complex $\Delta$ with face poset $P$,
let $M$ be a Morse matching of $\overline{P}$.
For~$i \geq 0$ 
let $u_i$ denote the number of critical $i$-dimensional
simplices.
Then
$\Delta$ is homotopy equivalent
to a CW-complex consisting of
$u_i$ $i$-cells, where $i \geq 0$. 
\end{theorem}

As it will be convenient for us to include the empty face in our
matching, we will use the following corollary of Forman's theorem.

\begin{corollary}
\label{corollary_wedge}
Let $\Delta$ be a simplicial
complex having Morse matching on the face
poset $P$ with exactly
$m$ critical $k$-dimensional simplices.
Then
$\Delta$ is homotopy equivalent
to a wedge
of $m$ $k$-dimensional spheres.
\end{corollary}

There is a straightforward criterion to determine
when a matching is a Morse matching.

\begin{lemma}
\label{lemma_general_Morse}
For a ranked poset $P$
to determine a given matching $M$ is a Morse matching,
it is enough to verify the Morse condition on
the
elements that are in the matching $M$ from adjacent
ranks in $P$.
\end{lemma}
\begin{proof}
We first show it is enough to check the
Morse acyclicity condition on elements from 
ranks $i$ and $i+1$ in the
poset $P$
for $i$ fixed.
As it is impossible to 
construct a cycle composed entirely
of edges oriented
downward,
any cycle 
in $P$ must include at least two elements from
the matching $M$.
So suppose $x$ and $y$ are two elements from the matching
of ranks $i$ and $i+1$ 
with the edge oriented from $x$ to $y$.
If the edge $x \rightarrow y$
is in a cycle, 
then the edge
$y \rightarrow z$
following this one
cannot
have
$\rho(z) = i+1$
(and hence, pointing upwards)
since
the 
element 
$y$ is already matched.
Similarly, if the edge $x \rightarrow y$
is in a cycle, 
then the edge 
$w \rightarrow x$ proceeding this one
cannot have $\rho(w) = i-1$ since
$x$ is already matched.
Hence we have shown if there is a cycle
in $P$, all the elements are from adjacent ranks.
Finally, if we have a cycle in $P$ on adjacent ranks,
the edges in the cycle alternate pointing upwards and downwards,
implying every other edge is from
the matching $M$ and hence all the elements in the
cycle are elements of $M$.
\end{proof}

What follows is a result which will be helpful when we
construct a Morse matching of the Whitehouse complex.

\begin{lemma}
\label{lemma_order_ideals}
Let $P$ be a ranked poset which
is the disjoint union
of a lower order ideal $L$
and
an upper order ideal $U$.
If $M_1$ is a Morse matching of $L$
and $M_2$ is a Morse matching of $U$,
then
$M_1 \cup M_2$ is a Morse matching of $P$.
\end{lemma}
\begin{proof}
By Lemma~\ref{lemma_general_Morse} 
without loss of generality
we
may assume all the elements we consider here
are 
those elements from  $M_1 \cup M_2$
of ranks $i$ and $i+1$, where $i \geq 1$.
As there is no element of $U$ matched with
an element of $L$,
all the edges from  rank $i+1$ elements in
$U$ to rank $i$ elements of $L$ are
oriented downwards.
Furthermore, as $U$ is an upper order ideal, 
there is no edge between any rank $i+1$ element in $L$
to a rank $i$ element of $U$.
Hence, it is impossible to construct a cycle
on the elements of $M_1 \cup M_2$,
implying
the matching $M_1 \cup M_2$ is indeed a Morse matching
of $P$.
\end{proof}

We now describe a Morse matching on the face poset
of the Whitehouse complex.

\begin{proposition}
\label{proposition_ABCD_poset}
Let $Q$ be the face poset
given by the images of
the maps $A$, $B$, $C$ and $D$
applied to the Whitehouse complex $\Delta_n$.
Let $M$ be the matching described by
orienting the edges from $A(F)$ to $B(F)$
and
$C(F,S) $ to $D(F,S)$,
where $S \subseteq F$ and $F$ ranges over all
faces of the Whitehouse complex $\Delta_n$.
Then 
$M$ is a Morse matching on $Q$.
\end{proposition}
\begin{proof}
Assume on the contrary that we can find a cycle between
rank $i$ and $i+1$ elements of $Q$.
The edges oriented from a rank $i$ to rank $i+1$ elements
in this cycle are Morse matched edges.
In terms of the forest representation of a face,
such an edge corresponds to
the element
$n+1$ being moved up one level higher in the tree.
(Here we are thinking of the leaves in the forest
as being the lowest level.)
In terms of the maps, such an edge corresponds
to $C(F,S)$ and $D(F,S)$
for some face $F$ and subset $S \subseteq F$
in $\Delta_n$.
The final move to raise the element $n+1$ corresponds
to the maps
$A(F)$ and $B(F)$.  
However, the path we have created cannot be continued,
and more importantly, cannot be completed to form a cycle,
since the element $n+1$ cannot be moved any higher.
Hence we cannot construct a cycle on the elements of $Q$,
so the matching $M$ described is a Morse matching.
\end{proof}

\begin{corollary}
The complex $Q$ described in Proposition~\ref{proposition_ABCD_poset}
is contractible.
\end{corollary}

We are now ready to prove our main result.
For other proofs, see~\cite{Robinson_Whitehouse, Sundaram, Vogtmann}.

\begin{theorem}
\label{theorem_pre-WDVV_homotopy}
The Whitehouse complex $\Delta_n$
is homotopy equivalent
to a wedge of $(n-2)!$
spheres
of dimension
$n-4$.
\end{theorem}
\begin{proof}
We proceed by
induction on the dimension $n$.
For $n = 3$,
the complex
$\Delta_3$ consists solely of the empty set,
so there is the trivial empty Morse matching.
This complex is homotopy equivalent to one
$(-1)$-dimensional sphere.

Begin to construct a Morse matching on the face poset
$P$ of 
$\Delta_{n+1}$ 
by
first 
orienting the edges in the
face poset  
$Q = \im(A(\Delta_n))\: \dot\cup \: \im(B(\Delta_n)) \: \dot\cup \: 
\im(C(\Delta_n)) \: \dot\cup \: \im(D(\Delta_n))$
as described in Proposition~\ref{proposition_ABCD_poset}.
The remainder of the face poset of $\Delta_{n+1}$
is $\im(E(\Delta_n))$.
Recall the image of $E$ applied to
$\Delta_n$ is isomorphic to $(n-1)$ copies of $\Delta_n$.
By induction, we have a Morse matching in each of these
$(n-1)$ copies of $\Delta_n$.
Each copy of $\Delta_n$ is an upper order ideal in
the face poset $P$.
Hence Lemma~\ref{lemma_order_ideals} applies,
so we have a constructed a Morse matching on the face poset of
$\Delta_{n+1}$.
In each of the $(n-1)$ copies
of $\Delta_n$ there
are $(n-2)!$ critical elements.
Moreover,
all the critical elements
are facets of dimension $(n-3)$.
Thus,
by Corollary~\ref{corollary_wedge}
the complex
$\Delta_{n+1}$ is homotopy equivalent to
a wedge of 
$(n-1)!$ spheres each having dimension $(n-3)$.
\end{proof}

\section{The Cohen-Macaulay property}
\label{section_Cohen_Macaulayness}

Recall that a simplicial complex $\Delta$
is {\em Cohen-Macaulay}
if the associated Stanley-Reisner ring
$k[\Delta]$ is Cohen-Macaulay~\cite{Stanley_green}.
In combinatorial commutative algebra
Cohen-Macaulay complexes have some very
nice enumerative properties.
Reisner's Criterion~\cite{Reisner, Stanley_green} gives a 
characterization of
Cohen-Macaulay simplicial complexes in terms of
reduced homology
of their links.

\begin{theorem}
A pure simplicial complex $\Delta$ 
is Cohen-Macaulay 
if and only if
for all faces $F \in \Delta$ and for
all $i < \dim(\link (F))$
we have
$\widetilde{H}_i(\link (F); k) = 0.$
\end{theorem}

We have the following result about the topology of the
links of the faces in the Whitehouse complex.

\begin{theorem}
\label{theorem_link_wedge}
For $F$ a face of 
the Whitehouse complex $\Delta_n$,
$\link_{\Delta_n}(F) $ is a wedge of 
$(c-1)! \cdot  \prod_{T \in F}(c(T) - 1)!$
spheres of 
dimension $n - 4 - |F|$,
where
$c$ is the number of components in the 
forest representation of $F$
and 
$c(T)$ is the number of children the node $T$ has in
the forest representation of $F$.
\end{theorem}
\begin{proof_special}
Let
$(S^n)^{\wedge k}$
denote the wedge of $k$ $n$-dimensional
spheres.
The free join of an $n$-dimensional sphere
with an $m$-dimensional sphere satisfies
$S^n * S^m \cong S^{n+m+1}$.
Additionally,
the wedge and free join operations
are distributive over
simplicial complexes,
that is,
$(X \wedge Y) * Z \cong (X * Z) \wedge (Y * Z)$.
For a proof of this fact, see~\cite[Lemma 3.14]{Dong}.
Hence
it follows that
$(S^n)^{\wedge k} * (S^m)^{\wedge l} \cong (S^{n+m+1})^{\wedge {k \cdot l}}$.
See also~\cite[Lemma 2.5 (ii)]{Bjorner_Welker}.
By Proposition~\ref{proposition_link}
and 
Theorem~\ref{theorem_pre-WDVV_homotopy}
we have
\begin{eqnarray*}
\link_{\Delta_n}(F) &\cong&   (S^{c-3})^{\wedge (c-1)!}
                             * \prod_{T \in F} 
                             (S^{c(T)-3})^{\wedge (c(T)-1)!}\\
                   &\cong&   \left(S^{|F| + c-3 + \sum_{T \in F} (c(T) - 3)}
                             \right)
                                   ^{\bigwedge (c-1)!  
                                     \cdot
                                     \prod_{T \in F}(c(T) - 1)!}\\
                   &\cong&   \left(S^{-2 |F| -3 + c +  \sum_{T \in F} c(T)}
                             \right)
                                   ^{\bigwedge (c-1)!  
                                     \cdot
                                     \prod_{T \in F}(c(T) - 1)!} \:\: .
\end{eqnarray*}
But
$c + \sum_{T \in F} c(T) = n-1 + |F|$.
Hence
\begin{eqnarray*}
                                        \:\:\:\: \:\:\:
                                        \:\:\:\: \:\:\:
                                        \:\:\:\: \:\:\:
                                        \:\:\:\: \:\:\:
                                        \:\:\:\: \:\:\:
                                         \:\:\:\: \:\:\:
                                        \:\:\:\: 
\link_{\Delta_n}(F) 
                   &\cong&   \left(S^{n-4-|F|} \right)
                                   ^{\bigwedge (c-1)!  
                                     \cdot
                                     \prod_{T \in F}(c(T) - 1)!} \:\: .
                                        \:\:\:\: \:\:\:
                                        \:\:\:\: \:\:\:
                                        \:\:\:\: \:\:\:
                                         \:\:\:\: \:\:\:
                                        \:\:\:\: \:\:\:
                                        \:\:\:\: \:\:\:
                                        \:\:\:\: 
                                        \qed
\end{eqnarray*}
\end{proof_special}

  From Reisner's Criterion and Theorem~\ref{theorem_link_wedge}
we have the following immediate result.
This can also be found in work of
Robinson-Whitehouse,
Sundaram and
Vogtmann~\cite{Robinson_Whitehouse, Sundaram, Vogtmann}.

\begin{theorem}
\label{theorem_Delta_n_is_Cohen_Macaulay}
The Whitehouse complex
$\Delta_n$
is Cohen-Macaulay
and hence the pre-WDVV ring
$\preWDVV_n$ is Cohen-Macaulay.
\end{theorem}

Trappmann and Ziegler~\cite{Trappmann_Ziegler}
 prove shellability of
the $k$-tree complex.  This is Hanlon's generalization of the Whitehouse tree
complex
corresponding to the case
$k=2$.
In unpublished work,
Wachs 
independently determined shellability of the Whitehouse
complex using a different shelling order.
As shellability implies
Cohen-Macaulayness, this gives another proof of 
Theorem~\ref{theorem_Delta_n_is_Cohen_Macaulay}.

\section{Concluding remarks}

%
%

\begin{table}
\begin{center}
\begin{tabular}{l|l}
$n$  &  $\Hilb(\WDVV_n)$\\  \hline
$3$ &   $1$\\
$4$ & $1 + t$\\
$5$ & $1 + 5t + t^2$\\
$6$ & $1 + 16t + 16 t^2 + t^3$\\
$7$ & $1 + 42t + 127 t^2 + 42t^3 + t^4$\\
$8$ & $1 + 99t + 715 t^2 + 715 t^3 + 99 t^4 + t^5$
\end{tabular}
\end{center}
\caption{The Hilbert series of the WDVV ring $\WDVV_n$ for
$3 \leq n \leq 8$.}
\label{table_WDVV_ring}
\end{table}

The author is currently studying the Hilbert series
of the WDVV ring.  The first few values are given 
in Table~\ref{table_WDVV_ring}.
The symmetry follows immediately from 
Poincar\'e duality and the fact the moduli
space $\overline{M_{g,n}}$
is smooth and compact.
It would be interesting to find
a natural combinatorial object corresponding to these values.

Vic Reiner has asked if the WDVV ring and the pre-WDVV
rings are Koszul.
Evidence for this is that both rings are defined by quadratic
relations and 
the reciprocal of the Hilbert series
for each ring
has alternating coefficients. From a result of Fr\"oberg~\cite{Froberg},
it follows that the Stanley-Reisner ring 
of a simplicial complex with minimal non-faces
having cardinality $2$ is Koszul.
This latter result applies to the pre-WDVV ring.

\begin{theorem}
The pre-WDVV ring is Koszul.
\end{theorem}

It remains to determine if either the pre-WDVV ring or
the WDVV ring are Gorenstein.

The complex of not $1$-connected graphs on 
$(n+1)$ vertices is also homotopy equivalent
to a wedge of $(n-2)!$ spheres
of dimension $n-4$; see~\cite{Babson_et_al, Vassiliev}.
Although the Whitehouse complex
$\Delta_n$ and this complex
are homotopy equivalent, 
they are not the same.
By Corollary~\ref{corollary_facets} the pre-WDVV complex is pure,
while the 
the cardinality of a facet 
in the complex of not $1$-connected graphs on
$(n+1)$ vertices
ranges from
$\lfloor \frac{n^2}{4}\rfloor$
to
${n \choose 2}$.

\section{Acknowledgements}

The author would like to thank 
Vic Reiner for pointing out references
and suggesting the Koszul question,
and
Louis Billera,
Richard Ehrenborg,
Arkady Kholodenko
and Michelle Wachs
for valuable comments and references.

Part of this research was done
while the author was a 
Visiting Member
at the Institute of Advanced Study
in Princeton during June~2002
and a Visiting Professor
at the University of Minnesota in August 2002.
This research was partially supported by 
a Summer 2002 University of Kentucky Faculty
Fellowship.

%
%
%

\newcommand{\journal}[6]{{\sc #1,} #2, {\it #3} {\bf #4} (#5), #6.}
\newcommand{\book}[4]{{\sc #1,} ``#2,'' #3, #4.}
\newcommand{\bookf}[5]{{\sc #1,} ``#2,'' #3, #4, #5.}
\newcommand{\thesis}[4]{{\sc #1,} ``#2,'' Doctoral dissertation, #3, #4.}
\newcommand{\springer}[4]{{\sc #1,} ``#2,'' Lecture Notes in Math.,
                          Vol.\ #3, Springer-Verlag, Berlin, #4.}
\newcommand{\preprint}[3]{{\sc #1,} #2, preprint #3.}
\newcommand{\preparation}[2]{{\sc #1,} #2, in preparation.}
\newcommand{\appear}[3]{{\sc #1,} #2, to appear in {\it #3}}
\newcommand{\submitted}[4]{{\sc #1,} #2, submitted to {\it #3}, #4.}
\newcommand{\DCG}{Discrete Comput.\ Geom.}
\newcommand{\JCTA}{J.\ Combin.\ Theory Ser.\ A}
\newcommand{\JCTB}{J.\ Combin.\ Theory Ser.\ B}
\newcommand{\AdvancesinMathematics}{Adv.\ Math.}
\newcommand{\JournalofAlgebraicCombinatorics}{J.\ Algebraic Combin.}
\newcommand{\communication}[1]{{\sc #1,} personal communication.}
\newcommand{\notes}[3]{{\sc #1,} Notes from Volume #2, dated
        #3.}
\newcommand{\collection}[9]{{\sc #1,} #2, 
           in {\it #3} (#4), #5,
           {\it #6}, {\bf #7}, #8, #9.}

\vspace{.5in}

{\sc
 \noindent
  Margaret\ A.\ Readdy                          \\
  Department of Mathematics                     \\
  University of Kentucky                        \\
  Lexington, KY 40506-0027                      \\
  {\rm readdy@ms.uky.edu}
}

\end{document}